\documentclass[a4paper, 12pt]{article}

\usepackage{graphicx}
\usepackage{url}
\usepackage{amsbsy,amssymb} 
\usepackage{amsmath}
\usepackage{abstract}
\usepackage{hyperref}
\usepackage{color}

\usepackage{amsmath,amsthm}

\usepackage[all]{xy} 

\newtheorem{Proposition}{Proposition}

\newtheorem{Theorem}{Theorem}
\newtheorem*{Proof}{Proof}
\newtheorem{Lemma}{Lemma}
\newtheorem{Definition}{Definition}

\begin{document}

{\LARGE\centering{\bf{The Poincare lemma for codifferential, anticoexact forms, and applications to physics}}}

\begin{center}
\sf{Rados\l aw Antoni Kycia$^{1,2,a}$}
\end{center}

\medskip
\small{
\centerline{$^{1}$Masaryk University}
\centerline{Department of Mathematics and Statistics}
\centerline{Kotl\'{a}\v{r}sk\'{a} 267/2, 611 37 Brno, The Czech Republic}
\centerline{\\}
\centerline{$^{2}$Cracow University of Technology}
\centerline{Faculty of Materials Engineering and Physics}
\centerline{Warszawska 24, Krak\'ow, 31-155, Poland}
\centerline{\\}

\centerline{$^{a}${\tt
kycia.radoslaw@gmail.com} }
}

\begin{abstract}
\noindent
The linear homotopy theory for codifferential operator on Riemannian manifolds is developed in analogy to a similar idea for exterior derivative. The main object is the cohomotopy operator, which singles out a module of anticoexact forms from the module of differential forms defined on a star-shaped open subset of a manifold. It is shown that there is a direct sum decomposition of a differential form into coexact and anticoexat parts. This decomposition gives a new way of solving exterior differential systems. The method is applied to equations of fundamental physics, including vacuum Dirac-K\"{a}hler equation, coupled Maxwell-Kalb-Ramond system of equations occurring in a bosonic string theory and its reduction to the Dirac equation.
\end{abstract}

\section{Introduction}

Well known formulation of the Poincar\'{e} lemma states \cite{SmoothManifoldsLee, ManifoldsTu} that 
\begin{equation}
H^{k}(\mathbb{R}^n) = H^{k}(point)=\left\{ \begin{array}{l}
                                       \mathbb{R}, \quad (k=0) \\
                                       0 \quad (k>0)
                                      \end{array}\right., 
\end{equation}
for $n>0$. This means that in $\mathbb{R}^n$ each differential form which is closed, i.e., it is in the kernel of exterior derivative operator $d$, is also exact, i.e., is also in the image of $d$. This can be also extended to a star-shaped open region of a smooth manifold $M$, which by definition is an open set $U$ of $M$ that is diffeomorphic to an open ball in $\mathbb{R}^n$, where $n=dim(M)$. In full generality, the Poincar\'{e} lemma is valid in any contractible submanifold.

In practical calculations, especially in physics, there is need for finding, if possible, potential from $d\alpha=0$ for some differential form $\alpha$, i.e., finding the exact formula for a differential form $\beta$ such that $d\beta=\alpha$. This can be done locally in a star-shaped region $U$ using (linear) homotopy operator \cite{EdelenExteriorCalculus, EdelenIsovectorMethods, ManifoldsTu, SmoothManifoldsLee, KyciaPoincare, DifferentialFormsInAlgebraicTopology}, i.e., 
 \begin{equation}
  H\omega := \int_{0}^{1} i_{\mathcal{K}} \omega|_{F(t,x)} t^{k-1} dt,
  \label{Eq.EdelenHomotopyOperator}
 \end{equation}
where a $k$-form $\omega \in \Lambda^{k}(U)$, $k=deg(\omega)$, 
\begin{equation}
 \mathcal{K}:=(x-x_{0})^{i}\partial_{i},
\end{equation}
and $F(t,x)=x_{0}+t(x-x_{0})$, $x_{0}\in U$, is a linear homotopy between the constant map $s_{x_{0}}:x\rightarrow x_{0}$ and the identity map $I:x\rightarrow x$. The form $\omega$ under the integral is evaluated at the point $F(t,x)$. Here $i_{\mathcal{K}}=\mathcal{K}\lrcorner $ is the insertion antiderivative. The assumption about star-shapedness of $U$ is imposed to ensure the correct definition of homotopy $F$ and the $H$ operator. Below we will be exclusively dealing with $U$ being star-shaped.

There is general construction of homotopy operator for various (non-linear) homotopies \cite{KyciaPoincare, ManifoldsTu, SmoothManifoldsLee, DifferentialFormsInAlgebraicTopology}, however the linear homotopy operator has additional advantage as it was pointed out by D.G.B. Edelen \cite{EdelenExteriorCalculus, EdelenIsovectorMethods, KyciaPoincare}.  The homotopy operator $H$ has many properties useful in calculations \cite{EdelenIsovectorMethods, EdelenExteriorCalculus}, e.g., 
\begin{equation}
 H^2=0, \quad HdH = H, \quad dHd =d,
 \label{Eq.H^2_HdH_dHd}
\end{equation}
\begin{equation}
 i_{\mathcal{K}}\circ H =0, \quad H\circ i_{\mathcal{K}}=0.
 \label{Eq.KH_HK}
 \end{equation}

This operator fulfils Homotopy Invariance Formula \cite{EdelenExteriorCalculus, EdelenIsovectorMethods, ManifoldsTu, SmoothManifoldsLee, KyciaPoincare}
\begin{equation}
 dH+Hd=I - s_{x_{0}}^{*},
 \label{Eq.EdelenHomotopyInvarianceFormula}
\end{equation}
where $s_{x_{0}}^{*}$ is the pullback along the constant map $s_{x_{0}}(x)=x_{0}$, and $I$ is the identity map.

It is well-known that the kernel of $d$ defines the closed vector space $\mathcal{E}(U)=\{\omega \in \Lambda(U) | d\omega=0 \}$ that is a vector subspace of $\Lambda(U)$. Since we focus on open star-shaped regions $U$ of a smooth manifold $M$ with or without boundary, by Poincar\'{e} lemma, elements of $\mathcal{E}(U)$ are also exact. Therefore, we will be using 'closed forms' and 'exact forms' interchangeably in what follows. 

Similarly, the kernel of $H$ on $U$ defines a module over $\Lambda^{0}(U)$ of antiexact forms $\mathcal{A}=\{\omega \in \Lambda(U)| H\omega=0\}$, which was described in \cite{EdelenIsovectorMethods, EdelenExteriorCalculus}. It was also proved \cite{EdelenIsovectorMethods, EdelenExteriorCalculus} that
\begin{equation}
\mathcal{A}=\{\omega \in \Lambda(U)|i_{\mathcal{K}} \omega=0, \quad \omega|_{x=x_{0}}=0\}, 
\end{equation}
and that there is a direct sum decomposition \cite{EdelenIsovectorMethods, EdelenExteriorCalculus, KyciaPoincare}
\begin{equation}
 \Lambda^{k}(U)=\mathcal{E}^{k}(U) \oplus \mathcal{A}^{k}(U),
\end{equation}
for $0 \leq k \leq n$.

On Riemannian manifolds with non-degenerate metric tensor $g$ one can define the Hodge star operator \cite{Nakahara, SpinorsBennTucker}, $\star:\Lambda^{r} \rightarrow \Lambda^{n-r}$,  that fulfils
\begin{equation}
 \star\star\omega = (-1)^{r(n-r)}\omega = (-1)^{r(n-r)} sig(g)\omega, \quad \omega \in \Lambda^{r}(U), 
\end{equation}
with the inverse
\begin{equation}
 \star^{-1} = (-1)^{r(n-r)} sig(g) \star = sig(g) \eta^{n-1}\star = sig(g)\star \eta^{n-1},
\end{equation}
where $\eta$ is an involutive automorophis: $\eta \omega = (-1)^{p}\omega$ for $\omega \in \Lambda^{p}$, and where $sig(g)=\frac{det(g)}{|det(g)|}$ is the signature of the metric $g$.
For clarity of presentation we will focus on the Riemannian case ($sig(g)=1$) only, and the other signatures, e.g., Lorentzian one, can be analyzed similarly.

Then the codifferential is defined as
\begin{equation}
 \delta = \star^{-1}d\star \eta.
\end{equation}

The Poincar\'{e} lemma for codifferential, an easy corollary of the Poincar\'{e} lemma, is as follows:
\begin{Theorem}(The Poincar\'{e} lemma for codifferential) \\
 For a star-shaped region $U$, if $\delta \omega=0$ for $\omega \in \Lambda^{k}(U)$, then there exists $\alpha \in \Lambda^{k+1}(U)$ for $k<n=dim(U)$, such that $\omega=\delta \alpha$.
\end{Theorem}

This paper aims to build a theory analogous to antiexact forms, in which codifferential $\delta$ is in the central place - anticoexact forms, and then apply it to various equations and systems of equations containing $d$ and $\delta$ operators. Anti(co)exact forms are local objects valid in a star-shaped open region of a manifold; however, local problems are essential to physics applications. Therefore, we also provide some examples of physics equations that can be solved locally by the methods presented here.

The paper is organized as follows: In the next section, we develop the theory of anticoexact forms that allows us to decompose arbitrary differential form into coexact and anticoexact parts. Then we relate this decomposition with exact-antiexact direct sum decomposition of \cite{EdelenExteriorCalculus, EdelenIsovectorMethods}. Next, the connection with Clifford algebras will be presented. Finally, application of (anti)(co)exact decomposition to various cases of Dirac(-K\"{a}hler) equations \cite{SpinorsBennTucker}, Maxwell equations of classical electrodynamics and their coupling with the Kalb-Ramond equations of bosonic string theory \cite{KalbRamondOriginalPaper, StringZwiebach} will be presented. In the Appendix the relation to the de Rham theory is discussed.

\section{Anticoexact forms}
This section defines an analog of the theory for antiexact forms, which we call anticoexact forms. The presentation will be along with Chapter 5 of \cite{EdelenExteriorCalculus} with marking differences between antiexact and defined below anticoexact forms.

We start from the homotopy operator for $\delta$:
\begin{Definition}
 We define the cohomotopy operator for $\delta$ for a star-shaped region $U$ as
 \begin{equation}
  h:\Lambda(U)\rightarrow \Lambda(U), \quad h=\eta \star^{-1}H\star.
  \label{Eq.hDefinition}
 \end{equation}
In particular,
 \begin{equation}
  h_{r}:\Lambda^{r}(U)\rightarrow \Lambda^{r-1}(U), \quad h_{r}=(-1)^{r+1}\star^{-1}H\star, \quad r>0.
 \end{equation}
Fig. \ref{Fig.ActionOfdHdeltah} presents interplay between all the operators.
 \begin{figure}
\centering
\xymatrix{ \ldots \ar@/_2pc/[r]_{d} \ar@/_/[r]_{h} & \ar@/_2pc/[l]_{H} \ar@/_/[l]_{\delta} \Lambda^{r-1}(U)  \ar@/_2pc/[r]_{d} \ar@/_/[r]_{h} & \ar@/_2pc/[l]_{H} \ar@/_/[l]_{\delta}  \Lambda^{r}(U) \ar@/_2pc/[r]_{d} \ar@/_/[r]_{h} & \ar@/_2pc/[l]_{H} \ar@/_/[l]_{\delta} \Lambda^{r+1} (U) \ar@/_2pc/[r]_{d} \ar@/_/[r]_{h} & \ar@/_2pc/[l]_{H} \ar@/_/[l]_{\delta} \ldots }
\caption{The action of $d$, $H$, $\delta$ and $h$ on $\Lambda$. Here $1<r<n-1$.}
\label{Fig.ActionOfdHdeltah}
\end{figure}
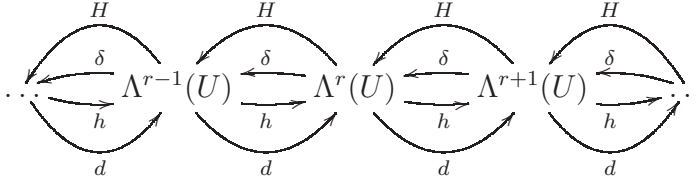
\end{Definition}
Since $(H\omega) |_{x=x_{0}}=0$, so $(h\omega)|_{x=x_{0}}=0$.

Such a definition makes the Homotopy Invariance Formula for $\delta$ and $h$ similar to (\ref{Eq.EdelenHomotopyInvarianceFormula}), namely,
\begin{Proposition}
 \begin{equation}
  \delta h + h \delta = I -S_{x_{0}},
  \label{Eq.CohomotopyInvarianceFormula}
 \end{equation}
 where $S_{x_{0}} =\star^{-1}s_{x_{0}}^{*}\star$. The operator $S_{x_{0}}$ is nonzero for $\Lambda^{n}(U)$ and it evaluates top forms at $x=x_0$, i.e., $S_{x_{0}}\omega = \omega|_{x=x_{0}}$ for $\omega \in \Lambda^{n}(U)$.
\end{Proposition}
\begin{Proof}
 Using the Homotopy Invariance Formula \ref{Eq.EdelenHomotopyInvarianceFormula} restricted to $\Lambda^{r}$, we have
 \begin{equation}
 \begin{array}{c}
  \star^{-1}(dH+Hd)\star = ((-1)^{r+1}\star^{-1}d\star) ((-1)^{r+1}\star^{-1}H\star) +  \\
  ((-1)^{r}\star^{-1}H\star)((-1)^{r}\star^{-1}d\star)= \delta h + h\delta = I - S_{x_{0}},
 \end{array} 
 \end{equation}
 since $\eta|_{\Lambda^{r}}= (-1)^{r} I$.
\end{Proof}

As a simple extension of the properties of $d$ and $H$ \cite{EdelenExteriorCalculus, EdelenIsovectorMethods}, we have
\begin{Proposition}
 \begin{equation}
  h^{2}=0, \quad \delta h \delta = \delta, \quad h \delta h = h.
 \end{equation}
\end{Proposition}
\begin{Proof}
 Since $H^{2}=0$ so, by (\ref{Eq.hDefinition}), $h^{2}=0$. For the second property, we have $\delta h \delta = \star^{-1}d\star\eta \eta \star^{-1}h \star \star^{-1}d*\eta = \star^{-1}d\star\eta =\delta$, since $\eta^{2}=1$ and $dHd=d$. Similarly, using $HdH=H$, we get the third property.
\end{Proof}

Define now the coclosed (that is also coexact in a star-shaped $U$) vector space
\begin{equation}
 \mathcal{C}:=\{\omega \in \Lambda(U) | \delta\omega =0\}.
\end{equation}
Note that $\mathcal{C}^{n} = \mathbb{R}\star 1$. Since $\mathcal{E}^{0}$ consists of constant functions \cite{KyciaPoincare}, coexact top forms are dual exact ones: $\star \mathcal{E}^{0}=\mathcal{C}^{n}$.

We have that 
\begin{Proposition}
 The operator $\delta h$ is the projector $\delta h: \Lambda \rightarrow \mathcal{C}$.
\end{Proposition}
\begin{Proof}
 For any form $\omega$, the form $\delta h \omega$ is coexat, so $\delta h:\Lambda \rightarrow \mathcal{C}$. It is idempotent since $\delta( h\delta h) = \delta h$. Finally, for $\omega\in \mathcal{C}$ we have from the Poincar\'{e} lemmat that there exists $\alpha$ such that $\delta \alpha =\omega$, and $h\omega = h\delta \alpha$. Then $\delta h \omega =\delta h \delta \alpha = \delta \alpha = \omega$. Therefore on $\mathcal{C}$ the operator $\delta h$ is the identity.
\end{Proof}

We can therefore define the coexact part of the form by the projection
\begin{equation}
 \omega_{c}:=\delta h\omega.
\end{equation}

By treating Homotopy Invariance Formula (\ref{Eq.CohomotopyInvarianceFormula}) as decomposition of the identity operator into projections, we can define 
\begin{Definition}(Anticoexact part of a form) \\
 We define anticoexact part of a form $\omega \in\Lambda^{k}(U), k<n$, on a star-shaped region $U$ as
\begin{equation}
 \omega_{ac}: = h\delta \omega =\omega - \delta h \omega.
\end{equation}
For $k=n$ the anticoexact part is 
\begin{equation}
 \omega_{ac}:=h\delta \omega = \omega - \omega|_{x=x_{0}}.
\end{equation}
\end{Definition}
Note that the anticoexact part of the form is of the type $\omega_{ac}=h\alpha$ and so these parts are in the kernel of the operator $h$ by its nilpotency.

We can define the anticoexact vector space as
\begin{equation}
 \mathcal{Y}(U):=\{ \omega \in \Lambda(U) | \omega = h\delta\omega \},
 \label{Eq.DefinitionY}
\end{equation}
that is the vector space of spanned by all anticoexact parts. Note that $\mathcal{Y}^{0} =0$.

Anticoexact space can be alternatively defined by the vector $\mathcal{K}$. To this end we have to use the following lemma
\begin{Lemma}(Equation (1.4.7) of \cite{SpinorsBennTucker})\\
 \begin{equation}
  i_{\alpha^{\sharp}}\star\phi = \star(\phi \wedge \alpha),
  \label{Eq.InsertionVsStar}
 \end{equation}
for $\alpha \in \Lambda^{1}$, $\phi$ an arbitrary form, and where $\sharp$ is a musical isomorphism such that $g(\alpha^{\sharp},X) = \alpha(X)$ for an arbitrary vector field $X$.
\end{Lemma}
Using this Lemma, we have
\begin{Proposition}
 \begin{equation}
  (\mathcal{K}^{\flat} \wedge) \circ h = 0, \quad h \circ \mathcal{K}^{\flat} \wedge =0.
  \label{Eq.kWedge_h_h_kWedge}
 \end{equation}
\end{Proposition}
\begin{Proof}
 Since we have (\ref{Eq.KH_HK}), i.e., $i_{\mathcal{K}}\circ H=0$, therefore $i_{\mathcal{K}} \star\star^{-1} \circ H=0$. Using (\ref{Eq.InsertionVsStar}), we have 
 $(\star\mathcal{K}^{\flat}\wedge) \circ \star^{-1}H=0$. From the definition (\ref{Eq.hDefinition}) of $h$ and the fact that $\star$ is an isomorphism, we have the result.
 
 For the second identity, using (\ref{Eq.KH_HK}), i.e., $H\circ i_{\mathcal{K}}=0$, we get $H\circ i_{\mathcal{K}} \star =0$, so using (\ref{Eq.InsertionVsStar}) we get $H\star \circ \mathcal{K}^{\flat}\wedge \circ \eta =0$. Therefore, $h\circ \mathcal{K}^{\flat}\wedge =0$, as required.
\end{Proof}

Using this we can characterize the vector space of anticoexact forms in an alternative way,
\begin{Proposition} {\ } \\ 
On a star-shaped region $U$,
\begin{equation}
 \mathcal{Y}(U)= \{\omega \in \Lambda(U)| \mathcal{K}^{\flat}\wedge \omega =0, \quad \omega|_{x=x_{0}}=0 \}.
 \label{Eq.DefinitionY_Kwedge}
\end{equation} 
\end{Proposition}
\begin{Proof}
 If $\omega \in \mathcal{Y}$, i.e., $\omega = h\delta \omega$ then from (\ref{Eq.kWedge_h_h_kWedge}) we get $\mathcal{K}^{\flat}\wedge \omega =0$ and $\omega|_{x=x_{0}}=0$.
 
 In the opposite direction, let $\mathcal{K}^{\flat}\wedge\alpha = 0$ and take $\omega = \alpha +\delta\beta$. Since $\mathcal{K}|_{F(t, x)}=t\mathcal{K}_{x}$ and
 $(i_{\mathcal{K}}\star\alpha)|_{F(t, x)}=t \star\mathcal{K}^{\flat}|_{x}\wedge\eta \alpha|_{F(t,x)}=0$, so $h\alpha=0$. Then $h\omega=h\delta\beta$, so $\omega_{c}=\delta h\omega=\delta h \delta \beta = \delta \beta$. Therefore the remaining part $\omega-\omega_{c}=\omega_{ac}=\alpha$. Since $\alpha$ is the anticoexact part of $\omega$ so $\alpha = h\delta \omega$ and therefore $\alpha_{x=x_{0}}=0$.
\end{Proof}

Contrary to $\mathcal{C}$ being only a vector space, we have
\begin{Proposition}  {\ } \\ 
 $\mathcal{Y}$ is a $C^{\infty}$-module.
\end{Proposition}
\begin{Proof}
 The conditions (\ref{Eq.DefinitionY_Kwedge}) defining $\mathcal{Y}$: $\mathcal{K}^{\flat}\wedge \omega =0$ and $\omega|_{x=x_{0}}=0$ is preserved under wedge multiplicatin of two elements from $\mathcal{Y}$ and under $C^{\infty}$ multiplication.
\end{Proof}

Antiexact forms can be written as $i_{\mathcal{K}}\alpha$ for some $\alpha$. Likewise, we have
\begin{Proposition} \label{Prop.AnticoexactKwedge}
 If $\omega \in \mathcal{Y}^{r}(U)$ then there exists $\alpha \in \Lambda^{r-1}(U)$ such that 
 \begin{equation}
  \omega = \mathcal{K}^{\flat}\wedge \alpha.
 \end{equation}
\end{Proposition}
\begin{Proof}
 Since $h\delta$ is the projector onto $\mathcal{Y}$, so for $\omega$ there is $\beta$ such that $\omega = h\beta$. Since $h$ is linear, so we can focus on a simple form which has local expression $\beta = f(x) dx^{I}$ for some multiindex $I$. Then 
 \begin{equation}
  h\beta = \eta \star^{-1}H\star\beta = \eta \star^{-1} i_{\mathcal{K}} \star dx^{I} \int_{0}^{1} dt f(F(t,x)) t^{|I|-1}  = \mathcal{K}^{\flat} \wedge \alpha,
 \end{equation}
where
\begin{equation}
 \alpha =  \left(\int_{0}^{1} dt f(F(t, x)) t^{|I|-1} \right)  dx^{I}.
\end{equation}
\end{Proof}

The final point of this section is the following
\begin{Theorem}
 For a star-shaped $U$ there is the direct sum decomposition
 \begin{equation}
  \Lambda^{k}(U) = \mathcal{C}^{k}(U) \oplus \mathcal{Y}^{k}(U).
 \end{equation}
\end{Theorem}
\begin{Proof}
 For $0<k<n$ from the Homotopy Invariance Formula (\ref{Eq.CohomotopyInvarianceFormula}), we have $h\delta + \delta h =I$. Moreover we know that both summands are projection operators. Therefore there is the unique decomopsition $\omega = \omega_{ac}+\omega_{c}$ with $\omega_{ac}=h\delta \omega$ and $\omega_{c}=\delta h \omega$.
 
 For $k=0$ we get $h\delta =I$ (i.e., $\Lambda^{0}(U)=\mathcal{Y}^{0}(U)$) since then $\delta \omega =0$.
 
 For $k=n$ we have $h\omega =0$, so $h\delta = I-S_{x_{0}}$. Therefore the decomposition is $\omega_{ac} = h\delta \omega = (I-S_{x_{0}})\omega$ and $\omega_{c} = S_{x_{0}}\omega \in \mathbb{R}\star 1=\mathcal{C}^{n}(U)$.
 
 Finally, if $\omega \in \mathcal{Y} \cap \mathcal{C}$ then using projectors $h\delta$ and $\delta h$ we get $\omega_{ac}=0=\omega_{c}$ and so $\omega =0$. As a result, the decomposition is indeed a direct sum decomposition.
 \end{Proof}

\section{(Anti)coexact vs. (anti)exact forms}
The theory developed in the previous section can be related to the theory of antiexact and exact forms, namely,
\begin{Proposition}
In a star-shaped region
 \begin{itemize}
  \item {$\mathcal{E}^{k} = \star \mathcal{C}^{n-k}$.}
  \item {$\mathcal{A}^{k} = \star \mathcal{Y}^{n-k}$}
 \end{itemize}
\end{Proposition}
\begin{Proof}
 If $\omega \in \mathcal{E}$, then $d\omega =0$. Therefore, $\delta \star \omega =0$, and so $\star \omega \in \mathcal{C}$. In a result $\star \mathcal{E} \subset \mathcal{C}$. Similarly one can prove the opposite inclusion. Since $\star$ operator is an isomorphism we get the first point.
 
 For the second point, we note that from (\ref{Eq.InsertionVsStar}) we can relate in a unique way $\mathcal{A}$ with $\mathcal{Y}$.
\end{Proof}

The proposition shows that the theory of exact and antiexact forms is dual, with respect to the $\star$ isomorphism, to the theory of coexact and anticoexact forms developed above. However, all the work done above is not futile. The notion of (anti)coexact forms and operator $h$ is helpful in applications, as we will see below.

Summing up the above results and the results of \cite{KyciaPoincare}, we have
\begin{Theorem}
\label{Th_decomposition}
 On a star-shaped open region $U$ of a manifold $M$, we have two ways of decompose $\Lambda(U)$ into the direct sums:
 \begin{itemize}
  \item {$\Lambda(U) = \mathcal{E}(U) \oplus \mathcal{A}(U)$;}
  \item {$\Lambda(U) = \mathcal{C}(U) \oplus \mathcal{Y}(U)$;}
 \end{itemize}
\end{Theorem}
Since the projection operators for these different decompositions do not commute, we cannot simultaneously decompose an element from $\Lambda(U)$ into both ways. The order of the decomposition is essential. 

In what follows, we need to recall the standard definitions.

\begin{Definition}
We will call $\mathcal{H}(U):=\mathcal{E}(U) \cap \mathcal{C}(U) = \{ \omega \in \Lambda(U) | d\omega=0 = \delta \omega\}$ Hodge harmonic forms \cite{GoldbergCurvature, Warner}. By analogy, we will call elements of $\bar{\mathcal{H}}:=\mathcal{A}\cap \mathcal{Y}$ Hodge antiharmonic forms.
\end{Definition}

In contrast, the typical (Kodaira) harmonic forms \cite{GoldbergCurvature, Warner} are elements of the kernel of the Laplace-Beltrami operator $\triangle=-(\delta d + d\delta)=(d-\delta)^2$. However, on non-compact $U$, these two definitions of harmonic forms are not in general equivalent \cite{Warner, GoldbergCurvature, deRham}. Due to the lack of compactness of $U$, the link with de Rham's theory is obscured, as will be explained in the Appendix.

In the next section the connection of above results with Clifford algebras is presented.

\section{Relation to Clifford algebras}

For a Riemannian manifold $(M,g)$ the Clifford bundle is isomorphic to $\Lambda(TM)$ pointwise by defining the Clifford multiplication of a vector $v \in T_{x}M$ by $\psi \in \Lambda(T_{x}M)$ as
\begin{equation}
 v\psi: = v \wedge \psi + i_{v}\psi,
 \label{Eq.CliffordMultiplication}
\end{equation}
see e.g., \cite{CliffordDifferentialForms, TuckerFermionsWithoutSpinors, TuckerDiracExterior, SpinorsBennTucker, SpinorialApproach, IndexTheory, ManyFacesOfMaxwell}. The conversion between elements of a tangent and cotangent bundle, if needed, is made using the metric $g$.

Then for the unique metric-compatible torsion-free connection $\bigtriangledown$ we can define for an orthonormal co-frame $\{e^{a}\}_{a=1}^{n}$ 
\begin{equation}
 d := e^{a}\wedge \bigtriangledown_{e_{a}}, \quad \delta := - i_{e^{a}}\bigtriangledown_{e_{a}}.
\end{equation}
In these terms, the Dirac(-K\"{a}hler) operator \cite{SpinorsBennTucker} on a Clifford bundle is defined as
\begin{equation}
 D := e^{a}\bigtriangledown_{e_{a}} = d - \delta.
\end{equation}

Note that for the form $\omega$ Clifford multiplied by $\mathcal{K}$ is
\begin{equation}
 \mathcal{K}\omega = \mathcal{K}^{\flat}\wedge \omega + i_{\mathcal{K}}\omega,
\end{equation}
which is the decomposition of $\omega$ into anticoexact and antiexact parts, see Proposition \ref{Prop.AnticoexactKwedge}.

In order to understand the structure of the Dirac operator on a Clifford bundle and its relation to the Poincar\'{e} lemma, we have to split it into grading of the base of Fig. \ref{Fig.ActionOfdHdeltah}. To this end, introduce the base of grading $\Lambda = \Lambda^{0} \oplus \ldots \oplus \Lambda^{n}$, where a form $\omega$ is written as the vector
\begin{equation}
 \omega = \left[ 
 \begin{array}{c}
  \omega^{0} \\
  \vdots \\
  \omega^{n}
 \end{array}
\right],
\end{equation}
where $\omega^{i} \in \Lambda^{i}$. In this base the exterior derivative $d$ has the simpler form
\begin{equation}
 d = \left[ 
 \begin{array}{ccccc}
  0 & 0 & \ldots & 0 & 0 \\
  d_{0} & 0 & \ldots & 0 & 0 \\
  0 & d_{1} & \ldots & 0 & 0 \\
  0 & 0 & \ddots & 0 & 0 \\
  0 & 0 & \ldots & d_{n-1} & 0 \\
 \end{array}
\right].
\end{equation}
Since $d_{k}d_{k-1}=0$ therefore $d^{2} =0$. The operator is nilpotent due to combination of the matrix multiplication and its (operator) elements. Likewise, we have
\begin{equation}
 \delta =  \left[ 
 \begin{array}{ccccc}
  0 & \delta_{1} & 0 & \ldots & 0 \\
  0 & 0 & \delta_{2} &\ldots & 0 \\
  0 & 0 & \ldots & \ddots & 0 \\
  0 & 0 & \ldots & 0 & \delta_{n}\\
  0 & 0 & \ldots & 0 & 0 \\
 \end{array}
\right],
\end{equation}
where $\delta^{2}=0$ by $\delta_{k-1}\delta_{k}=0$. 
For a star-shaped region $U$ we can define analogously the homotopy operators from Fig. \ref{Fig.ActionOfdHdeltah},
\begin{equation}
 H =  \left[ 
 \begin{array}{ccccc}
  0 & H_{1} & 0 & \ldots & 0 \\
  0 & 0 & H_{2} &\ldots & 0 \\
  0 & 0 & \ldots & \ddots & 0 \\
  0 & 0 & \ldots & 0 & H_{n}\\
  0 & 0 & \ldots & 0 & 0 \\
 \end{array}
\right], \quad
 h = \left[ 
 \begin{array}{ccccc}
  0 & 0 & \ldots & 0 & 0 \\
  h_{0} & 0 & \ldots & 0 & 0 \\
  0 & h_{1} & \ldots & 0 & 0 \\
  0 & 0 & \ddots & 0 & 0 \\
  0 & 0 & \ldots & h_{n-1} & 0 \\
 \end{array}
\right].
\end{equation}
Then the Dirac operator is
\begin{equation}
 D = d-\delta = \left[ 
 \begin{array}{cccccc}
  0 & -\delta_{1} & 0 & \ldots & 0  \\
  d_{0} & 0 & -\delta_{2} & \ldots & 0 \\
  0 & \ddots & \ddots & \ddots & 0 \\
  0 & \ldots & \ddots  & 0 & -\delta_{n} \\
  0 & \ldots & \ldots &  d_{n-1} & 0\\
 \end{array}
\right],
\end{equation}
with the Laplace-Beltrami operator
\begin{equation}
  D^{2} = (d-\delta)^{2} = -\left[ 
 \begin{array}{ccccc}
  \delta_{1}d_{0} & 0 & \ldots & 0 & 0 \\
  0 & d_{0}\delta_{1}+\delta_{2}d_{1} & 0 & \ldots & 0 \\
  0 & 0 & \ddots & 0 & 0 \\
  0 & \ldots & 0 & d_{n-2}\delta_{n-1} + \delta_{n}d_{n-1} & 0 \\
  0 & 0 & \ldots & 0 & d_{n-1}\delta_{n} \\
 \end{array}
\right].
\end{equation}
One then see that the Dirac operator mixes different grades and involves at most three neighbour grades.

Similarly, on $U$ we can define the anti-Dirac operator by $\mathcal{\reflectbox{D}}:=h-H$ and the anti-Laplace-Beltrami operator $\mathcal{\reflectbox{D}}^{2} = (h-H)^2 = -(hH+Hh)$.

A word on spinors realized by the Clifford bundle based on exterior bundle using (\ref{Eq.CliffordMultiplication}) is in order.  When $M$ is parallelizable, i.e., there is a frame $\{e_{a}\}_{a=1}^{n}$ such that $\bigtriangledown_{X} e_{a} =0$ for all $X\in \Gamma(TM)$ then from this basis one can construct global idempotents and use them to project Clifford algebra to minimal left ideals obtaining spinor subbundle of a Clifford bundle \cite{SpinorsBennTucker, CliffordDifferentialForms}. However, the existence of this ideal subbundle is more restrictive ($M$ must be parallelizable) than for the existence of the spinor bundle \cite{SpinorialApproach} in which case the vanishing of the second Stiefel-Whitney class is needed. The realization of spinors by (exterior) Clifford bundle gives, so-called, amorphous spinor fields/Dirac-Kh\"{a}ler spinors that do not have tensorial transformation law due to mixture of different grades, in contrast to covariant Dirac spinors used in physics. Therefore, they are not directly applicable to describe fermions (see however \cite{DiracKahler} for possible applications to lattice QCD), yet Dirac-like versions of Maxwell equations can be managed by methods presented in this paper. We will therefore focus only on sections of Clifford bundles -- Clifford fields. Detailed discussion of this subject is provided in \cite{ManyFacesOfMaxwell} and \cite{CliffordDifferentialForms, TuckerFermionsWithoutSpinors}.

In the next section, we examine local solutions of the Dirac and other physics equations using the machinery of decomposition of arbitrary form into (anti)(co)exact components.

\section{Application to equations of physics}
In this section, the application of the above theory to some equations of physics will be presented. All considerations will be given in a star-shaped region $U$ of a manifold $M$ since then it is possible to use the machinery presented above. The restriction to star-shaped open regions of a manifold sacrifice generality, however, for many applications in physics is usually sufficient. Global solutions that restrict to the local ones are usually connected with topological conditions on the manifold and require proper 'sheafication' procedure \cite{sheafication, GoldbergCurvature}. Therefore we will limit ourselves to the local considerations only.

\subsection{The vacuum Dirac equation}
First, we start from the solutions of the vacuum Dirac-K\"{a}hler equation \cite{DiracKahler, SpinorialApproach, SpinorsBennTucker} on a Clifford bundle
\begin{equation}
\mathcal{D}\psi=0
\label{Eq.VacuumDirac}
\end{equation}
in a star-shaped region $U$. The simplest solution is an arbitrary Hodge harmonic form $\psi \in \mathcal{H}(U)$, that is $\psi \in ker(d) \cap \ker(\delta)$. Such forms are also the solution of the Laplace equation $D^{2}\psi = 0$. On a compact manifold, $\psi$ is also a harmonic form. These solutions can be seen as a 'gauge modes' that allows to shift other solutions, since they nullify both terms $d$ and $\delta$ independently.

More complicated solutions involve three neighbour spaces $\Lambda^{k-1}$, $\Lambda^{k}$ and $\Lambda^{k+1}$ with $0<k<n$. Take two forms $\alpha \in \Lambda^{k-1}$ and $\beta\in\Lambda^{k+1}$ and set $\psi=\alpha+\beta$. Then the vacuum Dirac equation, under splitting into grades, gives the system
\begin{equation}
\left\{
 \begin{array}{c}
  \delta \alpha = 0 \\
  d \alpha - \delta \beta =0 \\
  d \beta =0.
 \end{array}
\right. 
\label{Eq.threeDiracequationVacuum}
\end{equation}

We have the following theorem
\begin{Theorem}
 The solution of (\ref{Eq.threeDiracequationVacuum}) for $0<k<n$ is either (\it{gauge case})
 \begin{equation}
  \alpha \in \mathcal{E}^{k-1}\cap \mathcal{C}^{k-1}, \quad \beta \in \mathcal{E}^{k+1}\cap\mathcal{C}^{k+1},
 \end{equation}
or (when $\alpha\not\in \mathcal{E}^{k-1}$ and $\beta \not\in\mathcal{C}^{k+1}$ - non-gauge case) is for
\begin{equation}
 d\alpha \in \mathcal{E}^{k}\cap \mathcal{C}^{k}, \delta\beta \in \mathcal{E}^{k}\cap\mathcal{C}^{k}.
\end{equation}
\end{Theorem}
\begin{Proof}
 From the first and the last equations of (\ref{Eq.threeDiracequationVacuum}) we get $\alpha \in \mathcal{E}^{k-1}$ and $\beta \in \mathcal{C}^{k+1}$.
 
 Using decomposition $dH+Hd=I$ for $0<k<n$ from the second equation of (\ref{Eq.threeDiracequationVacuum}) we get
 $$d\alpha - (dH+Hd)\delta\beta=0,$$
 that gives
 $$ d(\alpha-H\delta\beta)-Hd\delta\beta=0.$$
 This is the decomposition of $0$ into $\mathcal{E}\oplus\mathcal{A}$, so each component must vanish (we can also project into components using $Hd$ and $dH$ operators). We therefore have
 \begin{equation}
 \left\{ 
 \begin{array}{c}
  d(\alpha-H\delta\beta)=0 \\
  Hd\delta\beta=0.
 \end{array}
 \right.
 \label{Eq.DiracThreeProofSystem1}
 \end{equation}
 
If $\delta\beta=0$, i.e., $\beta \in \mathcal{C}$ then $d\alpha=0$ so $\alpha \in \mathcal{E}$ and we have
$$ \beta \in \mathcal{E}^{k+1}\cap \mathcal{C}^{k+1} \Rightarrow \alpha \in \mathcal{E}^{k+1}\cap \mathcal{C}^{k+1}.$$
 
Assume therefore that $\beta \not\in \mathcal{C}$, i.e., $\delta\beta \neq 0$. Then the second equation of (\ref{Eq.DiracThreeProofSystem1}) gives that $\delta\beta \in ker(Hd)=\mathcal{E}$, so $\delta \beta \in  \mathcal{E}^{k}\cap \mathcal{C}^{k}$. From the first equation of (\ref{Eq.DiracThreeProofSystem1}) and previous considerations $d\alpha \in \mathcal{E}^{k}\cap\mathcal{C}^{k}$.

Likewise, using decomposition $\delta h + h \delta =I$ for $0<k<n$ we obtain in the same way that either
$$ \alpha \in \mathcal{E}^{k+1}\cap \mathcal{C}^{k+1} \Rightarrow \beta \in \mathcal{E}^{k+1}\cap \mathcal{C}^{k+1},$$
or if $\alpha \not \in \mathcal{E}^{k-1}$ we get
$$d\alpha \in \mathcal{E}^{k}\cap \mathcal{C}^{k}, \delta\beta \in \mathcal{E}^{k}\cap\mathcal{C}^{k},$$
as required. 
\end{Proof}

Since the vacuum Dirac-K\"{a}ler equation is linear, so we can apply the above theorem to solve (\ref{Eq.VacuumDirac}) for $\psi$ involving more than three grades.

A similar procedure can be applied to the solutions of the vacuum anti-Dirac equation ($\mathcal{\reflectbox{D}}\psi =0$). The gauge modes are Hodge antiharmonic forms $\psi \in \bar{\mathcal{H}}$. The non-gauge solutions are of the form $\psi = \alpha + \beta$, where $\alpha \in \mathcal{A}^{k-1}\setminus\mathcal{Y}^{k-1}$ and $\beta \in \mathcal{Y}^{k+1}\setminus\mathcal{A}^{k+1}$ with the condition $h\alpha = H\beta \in \bar{\mathcal{H}}^{k}\subset\Lambda^{k}$.

\subsection{Massive Dirac equation}
The massive Dirac equation is of the form
\begin{equation}
 D\psi +\psi =0.
\label{Eq.MassiveDirac}
\end{equation}

For $\psi \in \Lambda^{k}$ belonging only to the one grade we have $\psi=0$. Therefore $\psi$ must be a sum of more grades. We consider a case of two neighbour grades $\psi=\alpha+\beta$ for $\alpha \in \Lambda^{k}$ and $\beta \in \Lambda^{k+1}$ for $0<k<n$. Then (\ref{Eq.MassiveDirac}) decomposes into the coupled system
\begin{equation}
 \left\{ 
 \begin{array}{c}
  \delta \alpha =0 \\
  \delta \beta = \alpha \\
  d\alpha = - \beta \\
  d\beta = 0,
 \end{array}
\right.
\label{Eq.MassiveDiracTwoGrades}
\end{equation}
which decouples into two constraints and two Klein-Gordon-type equations\footnote{For Lorentzian metric that is usual case in this context. For Euclidean metric the semilinear Laplace equations are obtained.}
\begin{equation}
 \left\{ 
 \begin{array}{c}
  \alpha \in \mathcal{C} \\
  \beta \in \mathcal{E} \\
  \triangle \alpha = \alpha \\
  \triangle \beta = \beta. \\
 \end{array}
\right.
\end{equation}

The above theory allows to reformulate (\ref{Eq.MassiveDiracTwoGrades}) as a coupled system of integral equations with constraints
\begin{equation}
 \left\{ 
 \begin{array}{c}
  \alpha \in \mathcal{C} \\
  \beta \in \mathcal{E} \\
  \alpha = -H\beta + dv \\
  \beta = h\alpha+\delta w,
 \end{array}
\right.
\label{Eq.MassiveDiracTwoGradesIntegral}
\end{equation}
where $v\in\Lambda^{k-1}$ and $w\in\Lambda^{k+2}$ for $k+2<n+1$ and zero otherwise, are two arbitrary forms.

Demanding that $\psi$ consists of components of all grades give a system that decouples into independent Klein-Gordon wave equations for all grades.

\subsection{Massless Dirac equation with source}
In the next step we will consider the massless Dirac equation
\begin{equation}
 D\psi = B,
\end{equation}
where $B \in \Lambda^{k}$, $0<k<n$ is a source term. In this section we will be looking for the solutions $\psi=\alpha+\beta$, where $\alpha \in \Lambda^{k-1}$ and $\beta \in \Lambda^{k+1}$. The asymmetric-grades ($k-1$ or $k+1$) cases are presented in the following sections since they represent electric and magnetic charges in Maxwell equations.

The system to solve is
\begin{equation}
 \left\{ 
 \begin{array}{c}
  \delta \alpha =0 \\
  d\alpha-\delta\beta = B \\
  d\beta =0.
 \end{array}
\right.
\label{Eq.MasslesSourceDirac}
\end{equation}
The first and the second equations give $\alpha \in \mathcal{C}^{k-1}$ and $\beta \in \mathcal{E}^{k+1}$ and we are left with the second equation. We have
\begin{Theorem}
 In order to solve (\ref{Eq.MasslesSourceDirac}) we can either

 \textbf{1st approach:}
 \begin{enumerate}
  \item {Solve $\triangle \beta = dB$ for $\beta$;}
  \item {Calculate $\alpha = H(\delta\beta+B)+dv$, where $v \in \Lambda^{k-2}$ if $k-2\geq 0$ and is arbitrary $(k-2)$-form, or $v=0$ for $k-2<0$;}
  \item {The solution is $\psi=\alpha+\beta$;}
 \end{enumerate}

 \textbf{2nd approach:}
 \begin{enumerate}
  \item {Solve $\triangle \alpha = -\delta B$ for $\alpha$;}
  \item {Calculate $\beta = h(d\alpha-B)+\delta w$, where $w \in \Lambda^{k+2}$ if $k+2\leq n$ and is arbitrary $(k-2)$-form, or $w=0$ for $k+2>n$;}
  \item {The solution is $\psi=\alpha+\beta$;}
 \end{enumerate}

\end{Theorem}

\begin{Proof}
We will focus on the middle equation of (\ref{Eq.MasslesSourceDirac}) and use the decompositions from Theorem \ref{Th_decomposition}.

First we use $dH+Hd=I$ for $1<k<n$ to get
$$ d\alpha-(dH+Hd)\delta\beta=(Hd+dH)B,$$
that is
$$d(\alpha-H\delta\beta-HB)-Hd(\delta\beta+B)=0,$$
which is nontivial $\mathcal{E}\oplus\mathcal{A}$-decomposition of $0$, that gives
\begin{equation}
\left\{ 
\begin{array}{c}
 \alpha = H\delta\beta+HB+d\phi \\
 \delta \beta + B \in \mathcal{E} \Leftrightarrow d\delta\beta=-dB,
\end{array}
\right.
\end{equation}
where $\psi$ is arbitrary. Since $d\beta=0$ we get the first case.

Likewise, we can use the decomposition $h\delta+\delta h =I$ for $1<k<n$ to get
$$ (h\delta+\delta h)d\alpha - \delta\beta = (h\delta+\delta h)B.$$
This gives nontivial $\mathcal{C}\oplus\mathcal{Y}$-decomposition of $0$, so 
\begin{equation}
 \left\{ 
 \begin{array}{c}
  \beta = hd\alpha-hB+\delta\psi \\
  d\alpha-B \in \mathcal{C} \Leftrightarrow \delta d\alpha-B=0,
 \end{array}
\right.
\end{equation}
where $\psi$ is arbitrary. By $\delta \alpha=0$ we get the second case. 
\end{Proof}

\subsection{Maxwell equations}
As a preparation for describing the Kalb-Ramond equations in the next subsection, we provide application of the above theory to the solutions of the Maxwell system on Minkowski space $M$ \cite{SpinorsBennTucker}
\begin{equation}
 dF = 0, \quad \delta F = j,
 \label{Eq.MaxwellSystem}
\end{equation}
where $F \in \Lambda^{2}$ and the external current is $j\in \Lambda^{1}$. This current is conserved since $\delta j=0$.

The typical approach on a star-shaped region $U$ is to take a potential $A \in \mathcal{A}^{1} \subset \Lambda^{1}$ such that $dA = F$, since $F\in \mathcal{E}^{2}$. Then the gauge transform $A\rightarrow A + \chi$, where $\chi\in \mathcal{E}^{1}$, that is $\chi = df$ for some $f\in\Lambda^{0}$, does not change $F$. The second equation is $-\delta d A = (\triangle + d\delta) A = -j$. Since $A \in \mathcal{A}^{1}$, so we can further decompose it into  $\mathcal{C}^{1}$ and $\mathcal{Y}^{1}$. We can remove anticoexact part of $A$ by imposing the Lorentz gauge: $\delta A =0$, and we obtain the wave equation  (in Lorentzian case)
\begin{equation}
 \triangle A = -j
 \label{Eq.MaxwellWave}
\end{equation}
that can be solved by standard propagator methods. Then there is still a gauge freedom $A \rightarrow A + \phi$, where $\triangle \phi = 0$.

Another approach is to use the above developed theory. First, use $dF=0$, i.e., $F\in \mathcal{E}^{2}$, to select, as before, $A \in \mathcal{A}^{1}$ such that $dA=F$. Then the second equation is $\delta dA =j$. Since the current is conserved, so $j \in \mathcal{C}^{1}$, and therefore, $j=\delta h j$. We have, $\delta (dA-hj)=0$, or $\delta(F-hj)=0$. In the end, the solution is 
\begin{equation}
 F= \delta \alpha +  hj,
 \label{Eq_F}
\end{equation}
where $\alpha \in \mathcal{Y}^{3}$ is defined up to an element of $\mathcal{C}^{3}$. Since (\ref{Eq_F}) represents the decomposition of $F$ into coexact and anticoexact parts, so by Theorem \ref{Th_decomposition} the element  $\delta\alpha$ is unique, and so $\alpha$ is unique up to a coexact form.  From this solution one sees that $j$ can be changed by an element of $\mathcal{Y}^{1}$ without affecting $F$. The additional constraint is $dF=0 = d\delta \alpha + dhj$, i.e., 
\begin{equation}
 \delta \alpha + hj \in \mathcal{E}^{2}.
 \label{Eq_F_constraint}
\end{equation}
The constraint is equivalent to $\delta \alpha + hj \in ker(Hd)$. 

The equation (\ref{Eq_F}) with the constraint (\ref{Eq_F_constraint}) is simply a decomposition of an element of $\mathcal{E}^2$ into $\mathcal{C}^2$ and $\mathcal{Y}^2$ according to Theorem \ref{Th_decomposition}.

We can further elaborate the condition (\ref{Eq_F_constraint}). Since $dF=0=d\delta \alpha + dhj$, we get, $\triangle \alpha+\delta d \alpha-dhj=0$. Removing from $\alpha$ exact part by imposing $d\alpha=0$, we get 
\begin{equation}
 \triangle \alpha =dhj,
 \label{Eq.MaxwellCodualWave}
\end{equation}
analogously to (\ref{Eq.MaxwellWave}) with a potential $A$ replaced with the three-form $\alpha$.

Note that, in this approach the existence of a specific $A$ was not needed - only the fact that $F\in \mathcal{E}^{2}$ is sufficient. In this approach a co-potential $\alpha$ is more important and it has also gauge freedom. Moreover, the current also can be modified by an element from  $\mathcal{Y}^{1}$ without affecting $F$. We can also recover $A$. Since $\delta \alpha + hj \in \mathcal{E}^{2}$, so $dA=F=dH(\delta \alpha + hj)$, and therefore, 
\begin{equation}
 A=df+H(\delta \alpha + hj),
\end{equation}
where $f\in \Lambda^{0}$.

This approach is more straightforward than that presented in \cite{EdelenExteriorCalculus} (Chapter 9) since we have the complete theory of (anti)exact and (anti)coexact forms at our disposal.

The full picture is well visible when we rewrite the system (\ref{Eq.MaxwellSystem}) using the Dirac operator \cite{SpinorsBennTucker}
\begin{equation}
 DF = -j.
\end{equation}
Then we can split $F = \psi + \gamma$, where $\psi = \alpha + \beta$ such that $D\psi=0$ is the solution of the vacuum Dirac equation with $\alpha \in \mathcal{C}^{1}$, $\beta \in \mathcal{E}^{3}$, and $\gamma \in \mathcal{E}^{2}$ is the solution of the nonhomogenous Dirac equation $D\gamma =-j$, that is $d\gamma=0$ and $\delta \gamma=j$.  

Note that if it would be that $j \in \Lambda^{3}$ (hypothetical magnetic monople current) and the equations would be $dF=j$, $\delta F=0$, then the procedure is similar as above with the restriction $F \in \mathcal{C}^{2}$. Since now $dj=0$, so $j=dHj$ and the first equation is $d(F-Hj)=0$, which gives $F = d \alpha + Hj$ for $\alpha\in \Lambda^{3}$ with the additional constraint $\delta(d \alpha +Hj)=0$, i.e., $d \alpha+Hj \in \mathcal{C}^{2}$. 
Since $\delta F=0$ so by the Poincar\'{e} lemma there exists $A\in \Lambda^{3}$ such that $F=\delta A$. Then the solution for $A$ is $A = \delta \beta + h(d\alpha + Hj)$ for some $\beta \in \Lambda^{4}$. This is dual to the classical electrodynamics presented above and by Theorem \ref{Th_decomposition} it is a decomposition of an element from $\mathcal{C}^2$ into a sum of elements from $\mathcal{E}^2$ and $\mathcal{A}^2$.

\subsection{Kalb-Ramond equations}
The Kalb-Ramond equations \cite{KalbRamondOriginalPaper, StringZwiebach} were postulated for describing charged bosonic string and, unlike Electrodynamics is the theory of a two-form $F$, they are equations for a three-form. The literature on the subject is vast, including both physical variations, e.g., \cite{KR1, KR2}, and generalizations of an idea of using $p$-forms, e.g., \cite{KR1Generalization, KR2Generalization, KR3Generalization, KR4Generalization, KR5Generalization}.

The equations have the following form
\begin{equation}
 dK = 0,\quad \delta K = J,
 \label{Eq.KalbRamondEquations}
\end{equation}
where $K\in\Lambda^{3}$ and $J\in \Lambda^{2}$. This can be further generalized to $p$-form electrodynamics \cite{KR3Generalization}, but we restrict ourselves to this simple example, since extension to different cases is straightforward.

From the first equation we have that $K \in \mathcal{E}^{3}$ and therefore there is a Kalb-Ramond field $B \in \Lambda^{2}$ such that $dB=K$. $B$ is defined up to $\mathcal{E}^{2}$, therefore we can chose it as $B\in\mathcal{A}^{2}$. From the second equation $J \in \mathcal{C}^{2}$ and therefore $J=\delta h J$ and so $\delta(K-hJ)=0$. We get, analogously to the electrodynamic, that $K = \delta \beta + hJ$, where $\beta \in \Lambda^{4}$ is defined up to gauge $\mathcal{C}^{4}$. The constraint is $\delta \beta + hJ \in \mathcal{E}^{3}$. Since in Minkowski space $d\beta=0$, so the constraint is the wave equation $\triangle \beta = dhJ$.

The equations (\ref{Eq.KalbRamondEquations}) can be written in the Dirac form 
\begin{equation}
 DK = -J,
\end{equation}
with the solution $K=\psi + \gamma$, where $D\psi=0$ and $\delta \gamma = J$, as in the case of Maxwell equations. 

Finally, we can couple Kalb-Ramond field with the Maxwell equation. It is possible since $B \in \mathcal{A}^{2}$ up to closed forms, and $F\in \mathcal{E}^{2}$. Therefore $F$ is a gauge field for $B$. We define \cite{StringZwiebach}
\begin{equation}
 R=B+F,
\end{equation}
which is possible since there is unique decomposition $\Lambda^{2} = \mathcal{E}^{2} \oplus \mathcal{A}^{2}$. Then we have $\delta R = \delta F + \delta B = j$ and therefore $\delta B=0$ by Maxwell equations (\ref{Eq.MaxwellSystem}). So $B \in \mathcal{A}^{2}\cap \mathcal{C}^{2}$. Taking exterior derivative $dR = dB = K$, and then $\delta K =\delta dB =J$.

Using the Dirac operator, the Kalb-Ramond-Maxwell system can be written in the compact form
\begin{equation}
 DR = K-j, \quad DK = -J,
\end{equation}
and can be analyzed similarly to the Maxwell equations.

\subsection{Cohomotopic fermionic harmonic oscillator}
In \cite{KyciaPoincare} there was presented a homotopy analogy of a fermionic quantum harmonic oscillator and its relation to Bittner's calculus of abstract derivative and integral \cite{R.Bittner}. Since operators $\delta$, $h$ are analogous to $d$ and $H$, therefore, we can define the co-version of this equation. The hamiltonian operator is
\begin{equation}
 \bar{H}:=h\delta-\delta h,
\end{equation}
and the anticommutator relation is played by nilpotency of $\delta$, $h$ and the Homotopy Invariance Formula (\ref{Eq.CohomotopyInvarianceFormula}) for $\Lambda^{r}$, where $0<r<n$. Then the eigenvalue problem is
\begin{equation}
 \bar{H}\omega = \lambda \omega, \quad \omega \in \Lambda^{r}
\end{equation}
where $\lambda \in \mathbb{R}$ are eigenvalues. As in \cite{KyciaPoincare}, the eigenvalues are $\lambda=\pm 1$ and eignevectors are coexact and anticoexact forms.

\section{Conclusions}
The theory of anticoexact forms in analogy to antiexact forms of \cite{EdelenExteriorCalculus} was developed in full detail. The most useful is the relation to Clifford algebra that allows us to solve the vacuum Dirac equation using (anti)(co)exact decomposition. Finally, the application of this decomposition in solving Maxwell equations of classical electrodynamics and Kalb-Ramond equations of bosonic string theory was presented. Moreover, (anti)(co)exact decomposition allows tracing all ingredients of the solutions of these and other similar equations.

\section*{Acknowledgments}
RK would like to thank Ioannis Chrysikos for the explanation of the rudiments of Clifford Algebras and Josef Silhan for support.

This research was supported by the GACR grant GA19-06357S and Masaryk University grant MUNI/A/0885/2019. RK also thank the SyMat COST Action (CA18223) for partial support.

\appendix

\section{Relation to de Rham theory}

Here we provide some hints on how to relate the above ideas to the de Rham/Hodge theory \cite{Warner, GoldbergCurvature, Morita, deRham}. The Poincar\'{e} lemma is a local theory that requires a star-shaped open region, and the de Rham theory is a global one which requires in a basic formulation a compact Riemannian manifold. Therefore, the intersection of these two will be visible on a manifold compactifying a star-shaped submanifold $U$. Since $U$ is diffeomorphic to an open disc $D$, so compactification requires two steps: making the closure of the disc $\bar{D}$ and then quotient out by the relation that glues boundary giving in the end a compact manifold. This is the case for sphere or torus. The differential forms that survive such procedure must behave well under the closure and then gluing. These two topological operations impose requirements on the elements from the spaces of (co)exact, anti(co)exact forms on $U$.
 
We will focus on a sphere $S^{n}$ as a compact connected manifold suitable for de Rham theory from one side. It can be seen as the one-point compactification of an open disc $U$ (or $\mathbb{R}^{n-1}$) that is star-shaped.

For compact connected manifolds, like $S^{n}$, the harmonic functions are constants. This agrees with the fact that closed forms on $U$ are constant forms \cite{KyciaPoincare} which extends smoothly to $S^{n}$ after the one-point compactification. 

However, for the (co)exact and anti(co)exact forms on $U$ under compactification to $S^{n}$ the nontrivial topology starts to play the role, and there is no easy correspondence to the forms from the de Rham theory. The apparent necessary condition for extending a differential form with continuous coefficients from $U$ to $S^{n}$ is the same limit value on the whole $\partial \bar{U}$ which is collapsed to the point under compactification. 

Since there are many possible ways to compactify even $\mathbb{R}^{n}$ (see Chapter 1 of \cite{Melorse}) such subject deserves another article.




\end{document}